\documentstyle[10pt]{article}
\begin{document}

\baselineskip 15pt
\parindent=1em
\hsize=12.3 cm \textwidth=12.3 cm
\vsize=18.5 cm \textwidth=18.5 cm

\def\supt{{\rm supt}}
\def\dom{{\rm dom}}
\def\bfone{{\bf 1}}
\def\Gen{{\rm Gen}}
\def\rk{{\rm rk}}
\def\cali{{\cal I}}
\def\calj{{\cal J}}
\def\Spec{{\rm Spec}}

\title{Understanding preservation theorems}
\author{
Chaz Schlindwein \\
Department of Mathematics and Computing \\
Lander University \\
Greenwood, South Carolina 29649, USA\\
{\tt cschlind@lander.edu}}

\maketitle

\def\forces{\mathbin{\parallel\mkern-9mu-}}
\def\notforces{\,\nobreak\not\nobreak\!\nobreak\forces}

\def\restr{\,\hbox{\vrule height8pt width.4pt depth0pt
   \vrule height7.75pt width0.3pt depth-7.5pt\hskip-.2pt
   \vrule height7.5pt width0.3pt depth-7.25pt\hskip-.2pt
   \vrule height7.25pt width0.3pt depth-7pt\hskip-.2pt
   \vrule height7pt width0.3pt depth-6.75pt\hskip-.2pt
   \vrule height6.75pt width0.3pt depth-6.5pt\hskip-.2pt
   \vrule height6.5pt width0.3pt depth-6.25pt\hskip-.2pt
   \vrule height6.25pt width0.3pt depth-6pt\hskip-.2pt
   \vrule height6pt width0.3pt depth-5.75pt\hskip-.2pt
   \vrule height5.75pt width0.3pt depth-5.5pt\hskip-.2pt
   \vrule height5.5pt width0.3pt depth-5.25pt}\,}

   \centerline{{\bf Abstract}}

   In its final form, this expository paper will present as much of Chapter VI of Shelah's book {\bf Proper and Improper Forcing} as I can manage.  Currently it has the special case of Theorem 1.12 giving the preservation of ${}^\omega\omega$-bounding and the preservation of the Sacks property.  I see no impediments to expanding this to an exposition of the full theorem.

\eject

\section{Introduction}

\section{Preservation of properness}

The fact that properness is preserved under countable support iterations was proved by Shelah in 1978.  The proof of this fact 
 is the basis of all preservation theorems for countable support iterations.

\proclaim Theorem 2.1 (Proper Iteration Lemma, Shelah). Suppose $\langle P_\eta\,
\colon\allowbreak\eta\leq\kappa\rangle$
is a countable support forcing iteration based on\/
$\langle\dot Q_\eta\,\colon\allowbreak\eta<\kappa\rangle$ and
for every $\eta<\kappa$ we have that\/ ${\bf 1}\forces_{P_\eta}``\dot Q_\eta$ {\rm is
proper.''} Suppose also that\/ $\alpha<\kappa$ and\/
$\lambda$ is a sufficiently large regular cardinal and\/ $N$ is a
countable elementary submodel of\/ $H_\lambda$ and\/
$\{P_\kappa,\alpha\}\in N$ and\/
$p\in P_\alpha$ is\/ $N$-generic and\/ {\rm
$p\forces``q\in  P_{\alpha,\kappa}\cap
N[G_{P_\alpha}]$.''}
Then there is $r\in P_\kappa$ such that
$r$ is $N$-generic and $r\restr\alpha=p$ and\/ {\rm
$p\forces``r\restr[\alpha,\kappa)\leq
q$.''}

Proof.  The proof proceeds by induction, so suppose that
the Theorem holds for all iterations of length less than $\kappa$.  Fix $\lambda$ a sufficiently large regular cardinal, and fix $N$ a countable elementary substructure of $H_\lambda$ such that $P_\kappa\in N$ and fix also $\alpha\in\kappa\cap N$ and $p\in P_\alpha$ and
a $P_\alpha$-name $q$ such that $p$ is $N$-generic and
$p\forces``q\in P_{\alpha,\kappa}\cap N[G_{P_\alpha}]$.''  

Case 1.  $\kappa=\beta+1$ for some $\beta$.

Because $\beta\in N$  we may use the induction hypothesis to fix $p'\in P_\beta$ such that
$p'\restr\alpha=p$ and $p'$ is $N$-generic and $p\forces``p'\leq q\restr\beta$.''  We have that
$p'\forces``q(\beta)\in N[G_{P_\beta}]$.'' Take $r\in P_\kappa$ such that
$r\restr\beta=p'$ and $p'\forces``r(\beta)\leq q(\beta)$ and $r(\beta)$ is
$N[G_{P_\beta}]$-generic for $Q_\beta$.''  Then $r$ is $N$-generic and we are done with the successor case.

Case 2.  $\kappa$ is a limit ordinal.

Let $\beta={\rm sup}(\kappa\cap N)$, and fix $\langle\alpha_n\,\colon\allowbreak
n\in\omega\rangle$ an increasing sequence from $\kappa\cap N$ cofinal in $\beta$ such that
$\alpha_0=\alpha$.  Let $\langle\sigma_n\,\colon\allowbreak n\in\omega\rangle$ enumerate all the $P_\kappa$ names $\sigma\in N$ such that ${\bf 1}\forces``\sigma$ is an ordinal.''

Using the induction hypothesis, build a sequence $\langle \langle p_n,q_n\rangle\,\colon\allowbreak n\in\omega\rangle$ such that $p_0=p$ and $q_0=q$ and for each $n\in\omega$ we have all of the following:

(1) $p_n\in P_{\alpha_n}$ and $p_n$ is $N$-generic  and
$p_{n+1}\restr\alpha_n=p_n$.

(2)  $p_n\forces``q_n\in P_{\alpha_n,\kappa}\cap N[G_{P_{\alpha_n}}]$ and
if $n>0$ then $q_n\leq q_{n-1}\restr[\alpha_n,\kappa)$ and $q_n\forces`\sigma_{n-1}\in
N[G_{P_{\alpha_n}}]$.'\thinspace''

(3) $p_n\forces``p_{n+1}\restr[\alpha_n,\alpha_{n+1})\leq q_n\restr\alpha_{n+1}$.''

Define $r\in P_\kappa$ such that $(\forall n\in\omega)\allowbreak
(r\restr\alpha_n=p_n)$ and ${\rm supp}(r)\subseteq\beta$.  
To see that $r$ is $N$-generic, suppose that $\sigma\in N$ is a $P_\kappa$-name for an ordinal. Fix $n$ such that
$\sigma = \sigma_n$. 
Because $p_{n+1}$ is $N$-generic, we have 
$p_{n+1}\forces``{\rm supp}(q_{n+1})\subseteq\kappa\cap N[G_{P_{\alpha_{n+1}}}]=
\kappa\cap N$,'' whence it is clear that 
$p_{n+1}\forces``r\restr[\alpha_{n+1},\kappa)\leq 
q_{n+1}$.'' 
We have $p_{n+1}\forces``q_{n+1}\forces`\sigma\in Ord\cap N[G_{P_{\alpha_{n+1}}}]=Ord 
\cap N
$,'\thinspace'' where $Ord$ is the class of all ordinals. 
Thus $r\forces``\sigma\in N$.''
We conclude that $r$ is $N$-generic, and the Theorem is established.

\proclaim Corollary 2.2 (Fundamental Theorem of Proper Forcing, Shelah).
Suppose $\langle P_\eta\,
\colon\allowbreak\eta\leq\kappa\rangle$
is a countable support forcing iteration based on\/
$\langle Q_\eta\,\colon\allowbreak\eta<\kappa\rangle$ and
for every $\eta<\kappa$ we have that\/ ${\bf 1}\forces_{P_\eta}`` Q_\eta$ {\rm is
proper.''}  Then $P_\kappa$ is proper.

Proof: Take $\alpha = 0$ in the Proper Iteration Lemma.

\section{Preservation of proper plus $\omega^\omega$-bounding}

In this section we recount Shelah's proof of the preservation of ``proper plus $\omega^\omega$-bounding.''  This is a special case of Theorem 1.12 of [PIF].  Other treatments of this material are Goldstern [Tools] and Goldstern and Kellner [forthcoming].

\proclaim Lemma 3.1.  Suppose $\langle P_\eta\,\colon\allowbreak
\eta\leq\kappa\rangle$ is a countable support iteration based on $\langle Q_\eta\,\colon\allowbreak \eta<\kappa\rangle$ and each $Q_\eta$ is proper in
$V[G_{P_\eta}]$.
Suppose ${\rm cf}(\kappa)=\omega$ and
 $\langle\alpha_n\,\colon\allowbreak n\in\omega\rangle$ is an increasing sequence of ordinals cofinal in $\kappa$ with $\alpha_0=0$.
Suppose also that $f$ is a $P_\kappa$-name for an element
of ${}^\omega\omega$, and suppose $p\in P_\kappa$.
Then there are
$\langle p_n\,\colon\allowbreak n\in\omega\rangle$ and $\langle f_n\,\colon\allowbreak n\in\omega\rangle$  such that
$p_0\leq p$ and for every $n\in\omega$ we have that each of the following holds:

(1) For all $k\leq n$ we have ${\bf 1}\forces_{P_{\alpha_n}}``p_0\restr[\alpha_n,\kappa)
\forces`f(k)=f_n(k)$.'\thinspace''

(2) $f_n$ is a $P_{\alpha_n}$-name for an element of ${}^\omega\omega$.

(3) ${\bf 1}\forces_{P_{\alpha_n}}``p_0\restr[\alpha_n,\alpha_{n+1})\forces`f_n(k)=f_{n+1}(k)$ for every $k\leq n+1$.'\thinspace''

(4) $p_{n+1}\leq p_n$.

(5) Whenever $k\leq m<\omega$ we have ${\bf 1}\forces_{P_{\alpha_n}}``p_m\restr[
\alpha_n,\alpha_{n+1})\forces`f_n(k)=f_{n+1}(k)$.'\thinspace''

\medskip

Proof:  
Fix $\lambda$ a sufficiently large regular cardinal and let $N$ be a countable
elementary substructure of $H_\lambda$ containing $P_\kappa$ and
$\langle\alpha_n\,\colon\allowbreak n\in\omega\rangle$ and $f$ and $p$.
Build $\langle p'_n\,\colon\allowbreak n\in\omega\rangle$ and
 $\langle\sigma_n\,\colon\allowbreak n\in\omega\rangle$ such that
$p'_0=p$ and each of the following holds:

(1) $p'_{n+1}\restr\alpha_n=p'_n\restr\alpha_n$.

(2) ${\bf 1}\forces_{P_{\alpha_n}}``\sigma_n\in\omega$ and
$p'_{n+1}\restr[\alpha_n,\kappa)
\forces`f(n)=\sigma_n$.'\thinspace''

(3) $p'_n\restr\alpha_n\forces``p'_{n+1}\restr[\alpha_n,\kappa)\leq 
p'_n\restr[\alpha_n,\kappa)$.''

(4) $p'_{n}\restr\alpha_n$ is $N$-generic.

(5) $p'_{n}\restr\alpha_n\forces``p'_{n}\restr[\alpha_n,\kappa)\in 
P_{\alpha,\kappa}\cap N[G_{P_{\alpha_n}}]$.''

(6)  $p'_{n}\in P_\kappa$.

Notice that (6) does not follow from the fact that $p'_{n}\restr\alpha_n\in
P_{\alpha_n}$ and $p'_{n}\restr\alpha_n\forces``p'_{n}\restr[\alpha,\kappa)\in
P_{\alpha,\kappa}$,'' but it does follow from (4) and (5).

Let $q_0=\bigcup\{p'_n\restr\alpha_n\,\colon\allowbreak n\in\omega\}$.  

At this point we define 
$f_n(k)=\sigma_k$ for $k\leq n$.  We have yet to define $f_n(k)$ for $k>n$.  Notice 
that we cannot set $f_n(k)=\sigma_k$ for $k>n$ because in $V[G_{P_{\alpha_n}}]$ 
we have that $\sigma_k$ is not an integer, but only a name.

It is easy to see that the following three properties hold:

(1') For all $k\leq n$ we have ${\bf 1}\forces_{P_{\alpha_n}}``q_0\restr[\alpha_n,\kappa)
\forces`f(k)=f_n(k)$.'\thinspace''

(2') $f_n\restr (n+1)$ is a $P_{\alpha_n}$-name for an element of ${}^{n+1}\omega$.

Fix $\lambda'$ a suffciently large regular cardinal and $M$ a 
countable elementary substructure of
$H_{\lambda'}$ containing $N$ and $q_0$.

Build $\langle p^n_0\,\colon\allowbreak n\in\omega\rangle$
and $\langle\tau_n\,\colon\allowbreak n\in\omega\rangle$
 such that $p^0_0=q_0$ and each of the following holds:

(1) $p_0^{n+1}\restr\alpha_n=p^n_0\restr\alpha_n$

(2) $p_0^{n+1}\leq p_0^n$

(3)  ${\bf 1}\forces_{P_{\alpha_n}}``\tau_n\in\omega$ and $
p_0^{n+1}\restr[\alpha_n,\alpha_{n+1})\forces`\tau_n=f_{n+1}(n+1)$.'\thinspace''

(4) $p^n_0\restr\alpha_n$ is $M$-generic.

(5) $p^n_0\restr\alpha_n\forces``p^{n+1}_0\restr[\alpha_n,\kappa)\in P_{\alpha_n,\kappa}\cap
M[G_{P_{\alpha_n}}]$.''

Notice that (1), (4), and (5) imply that $p_0^{n+1}\in P_\kappa$; this is the reason the
structure $M$ is needed.

There is no difficulty in doing this.  At this point, we define
$f_n(n+1)=\tau_n$ for every $n$.

Let $p_0=\bigcup\{p^n_0\restr\alpha_n\,\colon\allowbreak n\in\omega\}$.  

At this point the following parts of the Lemma are exemplified:

(1) For all $k\leq n$ we have ${\bf 1}\forces_{P_{\alpha_n}}``p_0\restr[\alpha_n,\kappa)
\forces`f(k)=f_n(k)$.'\thinspace''

(2') $f_n\restr(n+2)$ is a $P_{\alpha_n}$-name for an element of ${}^{n+2}\omega$.

(3) ${\bf 1}\forces_{P_{\alpha_n}}``p_0\restr[\alpha_n,\alpha_{n+1})\forces`f_n(k)=f_{n+1}(k)$ for every $k\leq n+1$.'\thinspace''

Choose $\lambda^*$ a sufficiently large regular cardinal.
We build $\langle p_n\,\colon\allowbreak n\in\omega\rangle$ 
and $\langle M_n\,\colon\allowbreak n\in\omega\rangle$ by recursion on $n\in\omega$.
Let $M_0$ be a countable elementary substructure of $H_{\lambda^*}$ containing
$M$ and $p_0$. 

Fix $n$, and suppose $p_n$ and $M_n$ have been defined.

For each $i<n$
let $q^i_n$ and $\xi^i_{n+1}$ be chosen such that
${\bf 1}\forces_{P_{\alpha_i}}``\xi^i_{n+1}\in\omega$ and
$q^i_n\in P_{\alpha_i,\alpha_{i+1}}\cap M_n[G_{P_{\alpha_i}}]$ and
$q^i_n\leq p_n\restr[\alpha_i,\alpha_{i+1})$ and
$q^i_{n}\forces`f_{i+1}(n+1)=\xi^i_{n+1}$.'\thinspace''

Build $\langle r^i_n\,\colon\allowbreak i\leq n\rangle$ such that for each $i\leq n$ we have the following.

(1) $r^i_n\in P_{\alpha_i}$.

(2) $r^i_n\leq p_n\restr\alpha_i$.

(3) $r^i_n$ is $M_n$-generic.

(4) If $i<n$ then $r^{i+1}_n\restr\alpha_i=r^i_n$.

(5) If $i<n$ then $r^i_n\forces``r^{i+1}_n\restr[\alpha_i,\alpha_{i+1})\leq q^i_n$.''

Then $p_{n+1}\restr\alpha_{n}=\bigcup\{r^i_n\,\colon\allowbreak i\leq n\}$ and
$p_{n}\restr[\alpha_{n},\kappa)=p_n\restr[\alpha_{n+1},\kappa)$. Let $M_{n+1}$ be a countable 
elementary substructure of $H_{\lambda^*}$ containing $M_n$ and $p_{n+1}$.

This completes the recursive construction.

We set $f_i(k)=\xi^i_k$ whenever $i+1<k$.

This completes the proof of the Lemma.

In the following Lemma we use the notation $``\check y$'' for the canonical $P$-name for $y$ when $y$ is a set in the ground model.

\proclaim Lemma 3.2.  Suppose $x$ is a $P$-name and $p\in P$ and
$p\forces``x\in V$.'' Then there is $q\leq p$ and $y$ such that
$q\forces``x=\check y$.''

Proof.  Well-known.

\proclaim Lemma 3.3.  Suppose $\langle P_\eta\,\colon\allowbreak\eta\leq\kappa\rangle$
is a countable support forcing iteration and $\lambda$ is a sufficiently large regular cardinal
and $N$ is a countable elementary substructure of $H_\lambda$ and $P_\kappa\in N$ and
$\alpha\in\kappa\cap N$ and $r\in P_\kappa$ is $N$-generic.  
Then\/ {\rm $r\restr\alpha\forces``r\restr[\alpha,\kappa)$ is
$N[G_{P_\alpha}]$-generic.''}

Proof: Suppose towards a contradiction that $\sigma$ is a $P_\alpha$-name for
a $P_{\alpha,\kappa}$-name for an ordinal and $r_1\leq r$ and $r_1\restr\alpha\forces``$the name
$\sigma$ is in $N[G_{P_\alpha}]$ yet $r_1\restr[\alpha,\kappa)\forces`\sigma\notin
N[G_{P_\alpha}]$.'\thinspace''
We may choose $r_2\leq r_1\restr\alpha$ and $x\in N$ a $P_\alpha$-name such that
$r_2\forces``$the $P_{\alpha,\kappa}$-names $\sigma$ and $x$ are equal.'' 
Viewing $x$ as a $P_\kappa$-name for an ordinal
rather than as a $P_\alpha$-name for a $P_{\alpha,\kappa}$-name for an ordinal, 
we may take an ordinal
$\tau$ and
$r_3\leq r_1$ such that $r_3\restr\alpha\leq r_2$ and $r_3\forces``x=\tau$.''
Because $r_3$ is $N$-generic, we have that $\tau\in N$.  Hence $r_3\forces``\tau\in N$.''
Hence $r_3\forces``\sigma\in N\subseteq N[G_{P_\alpha}]$.'' This is the desired contradiction.

\proclaim Definition 3.4.  For $f$ and $g$ in ${}^\omega\omega$ we say
$f\leq g$ iff $(\forall n\in\omega)\allowbreak(f(n)\leq g(n))$.
We say that $P$ is ${}^\omega\omega$-bounding iff 
$V[G_P]\models``(\forall f\in{}^\omega\omega)\allowbreak(\exists g\in {}^\omega\omega\cap V)
\allowbreak(f\leq g)$.''

\proclaim Theorem 3.5.  Suppose $\langle P_\eta\,\colon\allowbreak\eta\leq\kappa\rangle$ is a countable support iteration based on $\langle Q_\eta\,\colon\allowbreak\eta<\kappa\rangle$ and suppose\/ {\rm $(\forall\eta<\kappa)\allowbreak({\bf 1}\forces_{P_\eta}``Q_\eta$ 
is proper and ${}^\omega\omega$-bounding'').}
Then whenever $\lambda$ is a sufficiently large regular cardinal and\/ $N$ is a countable 
elementary substructure of $H_\lambda$
and $\alpha<\kappa$ and $\{P_\kappa,\alpha\}\in N$ and $p\in P_\alpha$ and $p$ is $N$-generic and $q$ and $f$ are $P_\alpha$-names in $N$ 
and\/ {\rm ${\bf 1}\forces_{P_\alpha}``q\in P_{\alpha,\kappa}$ and 
$f$ is a $P_{\alpha,\kappa}$-name
and $q\forces_{P_{\alpha,\kappa}}`f\in{}^\omega\omega$,'\thinspace''}  
then there are $q'$ and a $P_\kappa$-name $h$ such that 
{\rm $p\forces``q'\leq q$ and
$q'\forces`h\in{}^\omega\omega\cap V[G_{P_{\alpha}}]$ and $f\leq h$.'\thinspace''}

Proof: The proof proceeds by induction on $\kappa$.  We assume that $\lambda$, 
$N$, $\alpha$, $p$, $q$, and $f$ are as in the hypothesis of the Theorem.

Case 1. $\kappa=\beta+1$. 

 Because ${\bf 1}\forces_{P_\beta}``Q_\beta$ is ${}^\omega\omega$-bounding,''
we may take $q'$ and $h'$ to be $P_\beta$-names such that
${\bf 1}\forces_{P_\beta}``q^\#\leq q(\beta)$ and
$q^\#\forces`h'\in{}^\omega\omega\cap V[G_{P_\beta}]$ and $f\leq h'$.'\thinspace''
By Lemma 3.2 applied in $V[G_{P_\beta}]$, we may take $q^*$ and $h^*$ to be  
$P_\beta$-names such that ${\bf 1}\forces``q^*\leq q^\#$ and
$h^*\in{}^\omega\omega$ and $q^*\forces`h^*=h'$.'\thinspace''
We may assume that the names $q^*$ and $h^*$ are in $N$.
By the induction hypothesis we may take
$\tilde q$ and $h$ a $P_\beta$-name such  that
$p\forces``\tilde q\leq q\restr\beta$ and 
$q'\forces`h\in{}^\omega\omega\cap V[G_{P_\alpha}]$ and
$h^*\leq h$.'\thinspace''  Define $q'$ such that
$p\forces``q'=(\tilde q, q^*)\in P_{\alpha,\kappa}$.''
Clearly $p\forces``q'\leq q$ and $q'\forces`f\leq h$.'\thinspace''

Case 2.  ${\rm cf}(\kappa)>\omega$.

Because no $\omega$-sequences of ordinals can be added at limit stages of uncountable cofinality, we may take $\beta$ and $f'$ and $q'$ to be $P_\alpha$-names in $N$ such that
${\bf 1}\forces``\alpha\leq\beta<\kappa$ and $f'$ is a
$P_{\alpha,\beta}$-name and $q'\leq q$ and
$q'\restr\beta\forces_{P_{\alpha,\beta}}`f'\in {}^\omega\omega$ and
$q'\restr[\beta,\kappa)\forces_{P_{\beta,\kappa}}``f'=f$.''\thinspace'\thinspace''

For every $\beta_0\in\kappa\cap N$ such that $\alpha\leq\beta_0$ 
let $q^*(\beta_0)$ and $h(\beta_0)$ be  $P_\alpha$-names in $N$ such that
${\bf 1}\forces``$if $\beta=\beta_0$ and there is some $q^*\leq q'\restr\beta$ and some 
$h$ such that $h$ is a $P_{\alpha,\beta}$-name and
$q^*\forces`h\in{}^\omega\omega\cap V[G_{P_\alpha}]$ and 
$f'\leq h$,' then $q^*(\beta_0)$ and $h(\beta_0)$ are witnesses thereto.''
Let $q^*$ and $h$ and $s$ be $P_\alpha$-names such that for every $\beta_0
\in \kappa\cap N$, if $\alpha\leq\beta_0$, then
${\bf 1}\forces``$if $\beta=\beta_0$ then $q^*=q^*(\beta_0)$ and $h=h(\beta_0)$ and $s\in P_{\alpha,\kappa}$ and $s\restr\beta=q^*$ and
$s\restr[\beta,\kappa)=q'\restr[\beta,\kappa)$.''

Claim 1: $p\forces``s\leq q$  and $s
\forces`h\in{}^\omega\omega\cap V[G_{P_\alpha}]$ and $f\leq h$.'\thinspace''

Proof:  Suppose $p'\leq p$. Take $p^*\leq p'$ and $\beta_0<\kappa$ such that
$p^*\forces``\beta_0=\beta$.''  Because the name $\beta$ is in $N$ and $p^*$ is $N$-generic,
we have that $\beta_0\in N$.  Notice by the induction hypothesis that
we have $p\forces``$there is some $q^\#\leq q'\restr\beta_0$ and some $P_{\alpha,\eta_0}$-name
$h^\#$ such that $q^\#\forces`h^\#\in{}^\omega\omega\cap V[G_{P_\alpha}]$ and
$f'\leq h^\#$.'\thinspace''
Hence $p^*\forces``q^*=q^*(\beta_0)\leq q'\restr\beta$ and $h=h(\beta_0)$ and
$q^*\forces`h\in{}^\omega\omega\cap V[G_{P_\alpha}]$ and $f'\leq h$ and
$q'\restr[\beta,\kappa)\forces``f'=f$.''\thinspace'\thinspace''  Therefore
$p^*\forces``s\forces`f\leq h$.'\thinspace''

Claim 1 is established; this completes Case 2.

Case 3. ${\rm cf}(\kappa)=\omega$.

Let $\langle\alpha_n\,\colon\allowbreak n\in\omega\rangle$ be an increasing 
sequence from $\kappa\cap N$ cofinal in $\kappa$ such that 
$\alpha_0=\alpha$.

Let $\langle g_j\,\colon\allowbreak j<\omega\rangle$ list every $P_\alpha$-name
$g\in N$ such that
${\bf 1}\forces_{P_\alpha}`` g\in{}^\omega\omega$.''

Fix $\langle (p_n,f_n)\,\colon\allowbreak n\in\omega\rangle$ as in
Lemma 3.1 (applied in $V[G_{P_\alpha}])$.  That is, ${\bf 1}\forces``p_0\leq q$'' and 
 for every $n\in\omega$ we have that each of the following holds:

(0) $p_n$ is a $P_\alpha$-name for an element of $P_{\alpha,\kappa}$.

(1) For every $k\leq n$ we have ${\bf 1}\forces_{P_{\alpha_n}}``p_0\restr[\alpha_n,\kappa)
\forces`f(k)=f_n(k)$.'\thinspace''

(2) $f_n$ is a $P_{\alpha_n}$-name for an element of ${}^\omega\omega$.

(3) ${\bf 1}\forces_{P_{\alpha_n}}``p_0\restr[\alpha_n,\alpha_{n+1})\forces`f_n(k)=f_{n+1}(k)$ for every $k\leq n+1$.'\thinspace''

(4) ${\bf 1}\forces_{P_\alpha}``p_{n+1}\leq p_n$.''

(5) Whenever $k\leq m<\omega$ we have ${\bf 1}\forces_{P_{\alpha_n}}``p_m\restr[
\alpha_n,\alpha_{n+1})\forces`f_n(k)=f_{n+1}(k)$.'\thinspace''

We may assume that for every $n\in\omega$ the $P_\alpha$-names $p_n$ and $f_n$
are in $N$, and, furthermore, the sequence $\langle\langle 
p_n,f_n\rangle\,\colon\allowbreak n\in\omega\rangle$ is in $ N$.

In $V[G_{P_\alpha}]$, define $\langle g^n\,\colon\allowbreak n\in \omega\rangle$ by
$g^n(k)={\rm max}\{f_0(k),\allowbreak{\rm max}\{g_j(k)\,\colon\allowbreak
j\leq n\}\}$.
 Also in $V[G_{P_\alpha}]$ define $g\in{}^\omega\omega$ such that $g(k)=g^k(k)$ for all $k\in\omega$.

Claim 2.  Suppose $\alpha\leq\beta\leq\gamma<\kappa$ and $\{\beta,\gamma\}\in N$, and
suppose $f'\in N$ is a $P_\gamma$-name for an element of ${}^\omega\omega$. Suppose
$r\in P_\beta$ is $N$-generic and $r\forces``q'\in P_{\beta,\gamma}\cap N[G_{P_\alpha}]$.''
Then there are $q^*$ and $h$ such that $r\forces``q^*\leq q'$ and
$q^*\forces`h\in {}^\omega\omega\cap V[G_{P_\alpha}]$ and $f'\leq h$.'\thinspace''

Proof: In $V[G_{P_\alpha}]$ let $D=\{s\in P_{\alpha,\beta}\,\colon\allowbreak
s\forces``(\exists q^*\leq q')\allowbreak(\exists h\in{}^\omega\omega\cap V[G_{P_\alpha}])
\allowbreak(q^*\forces_{P_{\beta,\gamma}}`f'\leq h$')''$\}$.

Subclaim 1.  In $V[G_{P_\alpha}]$ we have that $D$ is a dense subset of
$P_{\alpha,\beta}$.

Proof: Given $\tilde p\in P_\alpha$ and $s$ a $P_\alpha$-name such that $\tilde
p\forces``s\in P_{\alpha,\beta}$,''
take $\lambda^*$ a sufficiently large regular cardinal and $M$ a countable elementary substructure of $\lambda^*$ containing $\tilde p$, $s$, and $N$.  Choose $p'\leq 
\tilde p$ such that $p'$ is $M$-generic.  By the overall induction hypothesis 
(i.e., because $\gamma<\kappa$), we have
$p'\forces``(\exists s'\in P_{\alpha,\gamma})\allowbreak(\exists h)\allowbreak
(s'\leq (s,q')$ and $h$ is a $P_{\alpha,\gamma}$-name and $s'\forces`h\in{}^\omega\omega
\cap V[G_{P_\alpha}]$ and $f'\leq h$').''  Fix such $s'$ and $h$.
By Lemma 3.2 we may take $h'$  and $\tilde s$ to be $P_\alpha$-names such that
$p'\forces``\tilde s\leq s'$ and 
$h'\in{}^\omega\omega$ and $\tilde s\forces`h'=h$.'\thinspace''
We have $p'\forces``\tilde s\restr[\alpha,\beta)\in D$.''
The Subclaim is established.

Working in $V[G_{P_\alpha}]$, let ${\cal J}\subseteq D$ be a maximal antichain 
of $P_{\alpha,\beta}$.
For each $s\in{\cal J}$ take $q^*(s)$ and $h(s)$ witnessing that $s\in D$.
Construct $q^\#$ and $h^\#$ such that $(\forall s\in{\cal J})\allowbreak
(s\forces_{P_{\alpha,\beta}}``q^\#=q^*(s)$ and $h^\#=h(s)$''). 

Subclaim 2.  $r\forces``q^\#\leq q'$ and $q^\#\forces_{P_{\beta,\gamma}}
`f'\leq h^\#$.'\thinspace''

Proof: Given $r_1\leq r$, take $r_2\leq r_1$ and a $P_\alpha$-name $s$ such that
$r_2\restr\alpha\forces``r_2\restr[\alpha,\beta)\leq s$ and $s\in{\cal J}$.''
We have $r_2\forces``q^\#=q^*(s)\leq q'$ and
$q^\#\forces`f'\leq h(s)=h^\#$.'\thinspace''

The Subclaim and the Claim are established.

Claim 3.  We may be build $\langle r_n\,\colon\allowbreak n\in\omega\rangle$ such that
$r_0=p$  and for every $n\in\omega$ we have that the following hold:

(1) $r_n\in P_{\alpha_n}$ is $N$-generic.

(2) $r_{n+1}\restr\alpha_n=r_n$.

(3) $r_n\forces``f_n\leq g$.''

Proof: Suppose we have $r_n$.

Take $P_{\alpha_n}$-names $F_0$ and $F_2$
such that ${\bf 1}\forces``$if there are functions $F_0'$ and $F_2'$
such that $F'_0$ maps $P_{\alpha_n,\alpha_{n+1}}$ 
into ${}^\omega\omega\cap V[G_{P_\alpha}]$ and
$F'_2$ maps $P_{\alpha_n,\alpha_{n+1}}$ into $P_{\alpha_n,\alpha_{n+1}}$ and
for every $q'\in P_{\alpha_n,\alpha_{n+1}}$ we have
$F'_2(q')\leq q'$ and $F'_2(q')\forces`f_{n+1}\leq F'_0(q')$', then $F_0$ and $F_2$ are witnesses
to this.''

We may assume that the names $F_0$ and $F_2$ are in $N$.

By Claim 2, we have $r_n\forces``F_0$ maps $P_{\alpha_n,\alpha_{n+1}}$ 
into ${}^\omega\omega\cap V[G_{P_\alpha}]$ and
$F_2$ maps $P_{\alpha_n,\alpha_{n+1}}$ into $P_{\alpha_n,\alpha_{n+1}}$ and
for every $q'\in P_{\alpha_n,\alpha_{n+1}}$ we have
$F_2(q')\leq q'$ and $F_2(q')\forces`f_{n+1}\leq F_0(q')$.'\thinspace''

In $V[G_{P_{\alpha_n}}]$, define $ g^*$ by
$g^*(i)= 
{\rm max}\{F_0(p_m\restr[\alpha_n,\alpha_{n+1}))(i)\,\colon\allowbreak m\leq i\}$.

We may assume the name $g^*$ is in $N$.

Notice that $r_n\forces``g^*\in N[G_{P_{\alpha_n}}]\cap V[G_{P_\alpha}] = N[G_{P_\alpha}]
$'' by Lemma 3.3.
Therefore we may choose a $P_{\alpha_n}$-name $k$ such that $r_n\forces``g^*=g_k$.''

Subclaim 1: $r_n\forces``F_2(p_k\restr[\alpha_n,\alpha_{n+1}))
\forces`f_{n+1}\leq g$.'\thinspace''

Proof: 
For $i\geq k$ we have 

\medskip

{\centerline{$r_n\forces``F_2(p_k\restr[\alpha_n,\alpha_{n+1}))
\forces`f_{n+1}(i)\leq F_0(p_k\restr[\alpha_n,\alpha_{n+1}))(i)
\leq g^*(i)=g_k(i)\leq g^i(i)=g(i)$.'\thinspace''}}

\medskip

\noindent The first inequality is by Subclaim 1, 
the second inequality is by the definition of $g^*$ 
along with the fact that $i\geq k$, the equality is by the 
definition of $k$, the next inequality is by the definition 
of $g^i$ along with the fact that $i\geq k$, and the last 
equality is by the definition of $g$.

For $i<k$, we have $r_n\forces``p_k\restr[\alpha_n,\alpha_{n+1})
\forces`f_{n+1}(i)=f_n(i)\leq g(i)$.'\thinspace''
The equality is by the choice of $\langle (f_m, p_m)\,\colon m\in\omega\rangle$ (see Lemma 3.1),
and the inequality is by the induction hypothesis that Claim 3 holds for integers less than or equal to $n$.

Because $r_n\forces``F_2(p_k\restr[\alpha_n,\alpha_{n+1}))\leq p_k\restr
[\alpha_n,\alpha_{n+1})$,'' we have that the Subclaim is established.

Using the Proper Iteration Lemma, take $r_{n+1}\in P_{\alpha_{n+1}}$ such that
$r_{n+1}$ is $N$-generic and $r_{n+1}\restr\alpha_n=r_n$ and
$r_n\forces``r_{n+1}\restr[\alpha_n,\alpha_{n+1})\leq F_2(p_k\restr[\alpha_n,\alpha_{n+1}))$.''

This completes the proof of Claim 3.

Let $r'=\bigcup\{r_n\,\colon\allowbreak n\in\omega\}$.  We have that
$p\forces``r'\restr[\alpha,\kappa)\leq q$ and 
$r'\restr[\alpha,\kappa)\forces`f\leq g$,'\thinspace''
 and so 
the Theorem is established.

\proclaim Corollary 3.6.  Suppose $\langle P_\eta\,\colon\allowbreak\eta\leq\kappa\rangle$ is a countable support iteration based on $\langle Q_\eta\,\colon\allowbreak\eta<\kappa\rangle$ and suppose\/ {\rm $(\forall\eta<\kappa)\allowbreak({\bf 1}\forces_{P_\eta}``Q_\eta$ 
is proper and ${}^\omega\omega$-bounding'').}  Then $P_\kappa$ is $\omega^\omega$-bounding.

Proof.  Take $\alpha = 0$ in Theorem 3.5.

\section{The Sacks property}

In this section we prove the preservation of ``proper plus Sacks property'' under countable support iteration.  The proof is due to Shelah; it is a special case of Theorem 1.12 of [PIF].  

\proclaim Definition 4.1.  
A poset $P$ has the Sacks property iff

{\centerline{${\bf 1}\forces_P``(\forall f\in{}^\omega\omega)\allowbreak(\exists H\in{}^\omega(\omega^{<\omega})\cap V)\allowbreak(\forall n\in\omega)\allowbreak(f(n)\in H(n))$.''}}

\proclaim Definition 4.2.  
For $T\subseteq{}^{<\omega}\omega$ a tree and $x\in{}^\omega(\omega-\{0\})$, 
we say that $T$ is an $x$-sized tree iff for every $n\in\omega$ we have that the 
cardinality of $T\cap{}^n\omega$ is at most $x(n)$.

\proclaim Definition 4.3.  For $T\subseteq{}^{<\omega}\omega$ we set $[T]$
equal to the set of all $f\in{}^\omega\omega$ such that
every initial segment of $f$ is in $T$.  That is, $[T]$ is the set of infinite branches of $T$.

\proclaim Definition 4.4.  For $x$ and $y$ in ${}^\omega(\omega-\{0\})$,
 we say that $x\ll y$ iff $(\forall n\in\omega)\allowbreak(x(n)\leq y(n))$ and
\[ \lim_{n\rightarrow\infty}y(n)/x(n)=\infty \]

\proclaim Lemma 4.5.  Suppose $y$ and $z$ are elements of ${}^\omega(\omega-\{0\})$ 
and $y\ll z$.  Suppose
that for each $n\in\omega$ we have 
 that\/ $T_n$ is a $y$-sized tree for every $n\in\omega$. 
 Then there is  a $z$-sized tree $T^*$ 
and an increasing sequence of integers $\langle k_n\,\colon\allowbreak n\in\omega\rangle$
such that $(\forall n\in\omega)\allowbreak(n<k_n)$ and
for every $f\in{}^\omega\omega$, we have

{\centerline{$(\forall n\in\omega)\allowbreak
(\exists i\leq n)\allowbreak(f\restr k_n\in T_{k_i})$ iff $f\in[T^*]$.}}

\medskip

Proof:   Build $\langle k_n\,\colon\allowbreak n\in\omega\rangle$ an increasing sequence of integers such that for every $n\in\omega$ we have $k_n>n$ and
$$(\forall t\geq k_{n})\allowbreak((n+2)y(t)\leq z(t))$$
Let $T^*=\{\eta\in{}^{<\omega}\omega\,\colon\allowbreak
(\forall n\in\omega)\allowbreak
(\exists i\leq n)\allowbreak(\eta\restr k_n\in T_{k_i})\}$.
Here, it is understood that if $m>{\rm length}(\eta)$ then $\eta\restr m$ denotes
$\eta$ itself.

It is easy to see that $T^*$ is nonempty and downward closed, and that every node in $T^*$ has at least one immediate successor in $T^*$. 
It is also easy to see that for all $f\in{}^\omega\omega$ we have
$(\forall n\in\omega)\allowbreak
(\exists i\leq n)\allowbreak(f\restr k_n\in T_{k_i})$ iff $f\in[T^*]$.

Claim. $T^*$ is a $z$-sized tree.

Proof.  Given $t\in\omega$, if $t<k_0$ then $T^*\cap {}^t\omega=T_{k_0}\cap {}^t\omega$, and hence $\vert T^*\cap{}^t\omega\vert\leq x_0(t)\leq z(t)$.  If instead
$t\geq k_0$, choose $n\in\omega$ such that $k_n\leq t<k_{n+1}$.
We have $ T^*\cap {}^t\omega=\{\eta\in{}^t\omega\,\colon\allowbreak(\forall j\leq n+1)\allowbreak(\exists i\leq j)\allowbreak(\eta\restr k_j\in T_{k_i})\}$ and so $\vert T^*\cap {}^t\omega\vert\leq\Sigma_{i\leq n+1}\vert T_{k_i}\cap {}^t\omega
\vert\leq(n+2)y(t)
\leq z(t)$.
The Claim and Lemma are established.

\proclaim Lemma 4.6.  Suppose  $y\in{}^\omega(\omega-\{0\})$
and $z\in{}^\omega(\omega-\{0\})$ 
and  $ y\ll z$.  
Suppose that for every $n\in\omega$ we have
that
$T_n$ is a $y$-sized tree and $h_n$ is an element of $[T_n]$. 
Suppose $f\in{}^\omega\omega$ and\/ {\rm $(\forall k\in\omega)\allowbreak
(f\restr k$ is an initial segment of $h_k)$.} 
Then there is a $z$-sized tree $T^*$ and a sequence of integers
$\langle m_i\,\colon\allowbreak i\in\omega\rangle$
such that
for every $\eta\in{}^{<\omega}\omega$, if $\eta$ is an initial segment of $f$
and $i\in\omega$ and
${\rm length}(\eta)\geq m_i$ and
 $\nu\in{}^{<\omega}\omega$ is an extension
of $\eta$ and $\nu\in T_{m_i}$, then  $\nu\in T^*$.

Proof.  By Lemma 4.5 we may choose $T^*$ a $z$-sized tree and
$\langle m_i\,\colon\allowbreak i\in\omega\rangle$ an increasing sequence of integers
such that

\medskip

{\centerline{$(\forall g\in{}^\omega\omega)\allowbreak((\forall n\in\omega)\allowbreak
(\exists i\leq n)\allowbreak(g\restr m_n\in T_{m_i})$ iff $g\in[T^*])$.}}

\medskip

Now suppose that $\eta\in{}^{<\omega}\omega$ and $i\in\omega$ and
${\rm length}(\eta)\geq m_i$ and
 $\eta$ is an
initial segment of $f$ and $\nu$ extends $\eta$
and $\nu\in T_{m_i}$.  We show $\nu\in T^*$.

Choose $h\in[T_{m_i}]$ extending $\nu$.  It suffices to show that $h\in[T^*]$.
Therefore it suffices to show that $(\forall k\in\omega)\allowbreak(\exists
j\leq k)\allowbreak(h\restr m_k\in T_{m_j})$.

Fix $k\in\omega$. We show that
$$h\restr m_k\in T_{{\rm min}(m_i,m_k)}$$
If $i\leq k$ then because $h\in[T_{m_i}]$ we have that
$h\restr m_k\in T_{m_i}$ and we are done.  If instead $k<i$ then
$h\restr m_k=\eta\restr m_k=f\restr m_k$ which is an initial segment of $h_{m_k}\in
[T_{m_k}]$.  Hence $h\restr m_k\in T_{m_k}$.

QED Lemma 4.6.

\proclaim Lemma 4.7.  Suppose $f\in{}^\omega\omega$ and 
$y\in{}^\omega(\omega-\{0\})$ and $z\in{}^\omega(\omega-\{0\})$, and suppose
$\langle x_n\,\colon\allowbreak n\in\omega\rangle$ is a sequence of elements
of ${}^\omega(\omega-\{0\})$ such that $(\forall n\in\omega)\allowbreak
(x_n\ll x_{n+1}\ll y\ll z)$.
Suppose for every $n\in\omega$, we have $x_n^*\in{}^\omega(\omega-\{0\})$ and 
$x_n^*\ll x_n\ll x^*_{n+1}$, and we have
$\langle x_{n,j}\,\colon\allowbreak n\in\omega$, $j\in\omega\rangle$ is a
sequence of elements of ${}^\omega(\omega-\{0\})$ such that
for every  $j\in\omega$ we have
$x_n\ll x_{n,j}\ll x_{n,j+1}\ll x_{n+1}^*$.
Suppose  $\langle T_{n,j}\,\colon\allowbreak n\in\omega$, $j\in\omega\rangle$ is
a sequence such that for every $n\in\omega$ and
 $j\in\omega$ we have that $T_{n,j}$ is an
$x_{n,j}$-sized tree.  Then there are $\langle T^n\,\colon\allowbreak n\in\omega\rangle$ and
$T^*$ such that $T^*$ is a $z$-sized tree and $f\in[T^0]$ and $T^0\subseteq T^*$ and for every
$n\in\omega$ we have

(i) $T^n\subseteq T^{n+1}$ and $T^n$ is an $x_n$-sized tree, and

(ii) for every $j\in\omega$ and every $g\in[T_{n,j}]$ 
there is $k\in\omega$ such that for every $\eta\in{}^{<\omega}\omega$ extending $g\restr k$, if 
$\eta\in T_{n,j}$ and $\eta\restr k\in T^n\cap T^*$ then
$\eta\in T^{n+1}\cap T^*$.

\medskip

Proof:    Choose $T^0$ an
$x_0$-sized tree such that $f\in[T^0]$.  
Given $T^n$, build $\langle T'_{n,j}\,\colon\allowbreak j\in\omega\rangle$
as follows.  Let $T'_{n,0}=T^n$.  Given
$T'_{k,j}$ take $m(n,j)\in\omega$ such that
$(\forall t\geq m(n,j))\allowbreak(2x_{n,j}(t)\leq x_{n,j+1}(t))$.
Let $T'_{n,j+1}=\{\eta\in T_{n,j}\colon\allowbreak
\eta\restr m(n,j)\in T'_{n,j}\}\cup T'_{n,j}$.

Claim 1.  Whenever $i\leq j<\omega$ we have $T'_{n,i}\subseteq T'_{n,j}$.

Proof.  Clear.

Claim 2. Suppose $T^n$ is an $x_n$-sized tree.  Then $(\forall j\in\omega)\allowbreak(T'_{n,j}$
is an $x_{n,j}$-sized tree).

Proof:  It is clear that $T'_{n,0}$ is an $x_{n,0}$-sized tree.  
Assume that $T'_{n,j}$ is an $x_{n,j}$-sized tree.  Fix $t\in\omega$.  
If $t<m(n,j)$ then $T'_{n,j+1}\cap {}^t\omega=
T'_{n,j}\cap{}^t\omega$ and so $\vert T'_{n,j+1}\cap {}^t\omega\vert\leq x_{n,j}(t)\leq
x_{n,j+1}(t)$.  So assume that $t\geq m(n,j)$.  
We have $T'_{n,j+1}\cap{}^t\omega\subseteq(T'_{n,j}\cap{}^t\omega)\cup (T_{n,j}\cap
{}^t\omega)$.  Therefore
we have $\vert T'_{n,j+1}\cap{}^t\omega\vert\leq 2x_{n,j}(t)\leq x_{n,j+1}(t)$.  The
Claim is established.

For each $n\in\omega$, using Claim 2 and Lemma 4.5 we my find an increasing sequence of integers $\langle k_{n,j}\,\colon\allowbreak j\in\omega\rangle$ and $T^{n+1}$ such that
$(\forall j\in\omega)\allowbreak(k_{n,j}>j)$ and
if $T^n$ is an $x_n$-sized tree, then $T^{n+1}$ is
an $x_{n+1}$-sized tree such that
for all $\eta\in{}^{<\omega}\omega$, we have 

\medskip

{\centerline{$(\forall j\in\omega)\allowbreak
(\exists i\leq j)\allowbreak
(\eta\restr k_{n,j}\in T'_{n,k_{n,i}})$ iff $\eta\in T^{n+1}.$}}

\medskip

Applying mathematical induction, we have that each $T^n$ is in fact an $x_n$-sized tree.

Claim 3. $T^n\subseteq T^{n+1}$ for every $n\in\omega$.

Proof: By Claim 1 we have that $T^n\subseteq T'_{n,i}$ for every $i\in\omega$.
By the definition of $T^{n+1}$ we have that
$T^{n+1}\supseteq\bigcap\{T'_{n,k_{n,i}}\,\colon\allowbreak i\in\omega\}\supseteq
\bigcap\{T'_{n,i}\,\colon\allowbreak i\in\omega\}\supseteq T^n$.

Applying Lemma 4.5 again we obtain an increasing sequence of integers $\langle k_n\,\colon\allowbreak n\in\omega\rangle$ and a
$z$-sized tree $T^*$ such that $(\forall n\in\omega)\allowbreak
(n<k_n)$ and for every $\eta\in{}^{<\omega}\omega$,
we have that

\medskip

{\centerline{$(\forall n\in\omega)(\exists i\leq n)(
\eta\restr k_n\in T^{k_i})$ iff $\eta\in T^*$.}}

\medskip

Notice that $T^0\subseteq\bigcap\{T^n\,\colon\allowbreak n\in\omega\}\subseteq T^*$.

Now we verify that $\langle T^n\,\colon\allowbreak n\in\omega\rangle$ and $T^*$ satisfy the
remaining conclusions of the Lemma.  Accordingly, fix $n\in\omega$ and $j\in\omega$ and
$g\in[T_{n,j}]$.
Let $k={\rm max}(k_n,\allowbreak{\rm max}\{k_{n,j'}\,\colon\allowbreak j'\leq j\},
\allowbreak{\rm max}\{m(n,j')\,\colon\allowbreak
j'\leq j\})$.  Fix $\eta\in{}^{<\omega}\omega$ extending $g\restr k$ and assume that
$\eta\in T_{n,j}$ and $\eta\restr k\in T^n\cap T^*$.

Claim 4. $\eta\in T^{n+1}$.

Proof:  It suffices to show $(\forall j'\in\omega)\allowbreak
(\exists i\leq j')\allowbreak(\eta\restr k_{n,j'}\in T'_{n,k_{n,i}})$.  Fix $j'\in\omega$
and let $i={\rm min}(j,j')$.

Case 1: $j'\leq j$. 

Because $k_{n,j'}\leq k$ we have that 
 $\eta\restr k_{n,j'}\in T^n\subseteq T'_{n,k_{n,i}}$, as required.

Case 2: $j<j'$.

It suffices to show that $\eta\restr k_{n,j'}\in T'_{n,k_{n,j}}$. 
Because $g\restr k\in T^n$ we have that
 $g\restr m(n,j)\in T^n\subseteq T'_{n,j}$.  Because
we have $\eta\in
T_{n,j}$ and $\eta\restr m(n,j)=g\restr m(n,j)\in T'_{n,j}$, we know
 by the definition of
$T'_{n,j+1}$ and Claim~1 that
$\eta\in T'_{n,j+1}\subseteq T'_{n,k_{n,j}}$.

Claim 5. $\eta\in T^*$.

Proof: It suffices to show $(\forall i\in\omega)\allowbreak(\exists i'\leq i)
\allowbreak(\eta\restr k_i\in T^{k_{i'}})$.  Towards this end, fix $i\in\omega$.

Case 1:  $i\leq n$.

Because $\eta\restr k\in T^*$ and $\eta$ extends $g\restr k$,
 we have $g\restr k_i\in T^*$ and hence
we may take $i'\leq i$ such that $g\restr k_i\in T^{k_{i'}}$. But we also have
 $\eta\restr k_i=g\restr k_i$, so we finish case~1.

Case 2: $n<i$.

We let $i'=i$. We have
$\eta\restr k_i=g\restr k_i\in T^n$.  Therefore
by Claim 4 we have $\eta\restr k_i\in T^{n+1}$, and by Claim 3 we have that
$T^{n+1}\subseteq T^{k_n}\subseteq T^{k_i}$.

The Lemma is established.

\proclaim Lemma 4.8.  Suppose $N$ is a countable set and $z\in{}^\omega(\omega-\{0\})$ and $z$ eventually dominates each element of
${}^\omega(\omega-\{0\})\cap N$.  Then there is $y\ll z$ such that $y$ eventually dominates each element of  ${}^\omega(\omega-\{0\})\cap N$. 

Proof.  Let $\langle x_n\,\colon\allowbreak n\in\omega\rangle$ list
${}^\omega(\omega-\{0\})\cap N$ such that $x_0(0)=1$. For each $n\in\omega$ let
$t_n$ be an integer such that for every $t\geq t_n$ we have
$z(t)\geq n^2\,x_n(t)$.  We may assume $t_0=0$.  For every $t\in\omega$ set
$y(t)={\rm max}\{ n x_n(t)\,\colon\allowbreak t_n\leq t\}$. 
It is easy to verify that this works.

\proclaim Theorem 4.9.  Suppose $\langle P_\eta\,\colon\allowbreak\eta\leq\kappa\rangle$ is a countable support iteration based on $\langle Q_\eta\,\colon\allowbreak\eta<\kappa\rangle$ and suppose\/ {\rm $(\forall\eta<\kappa)\allowbreak({\bf 1}\forces_{P_\eta}``Q_\eta$ 
is proper and has the Sacks property'').}
Then whenever $\lambda$ is a sufficiently large regular cardinal and\/ $N$ is a countable 
elementary substructure of $H_\lambda$
and $\alpha<\kappa$ 
and
$\{P_\kappa,\alpha\}\in N$ and $p\in P_\alpha$ and $p$ is $N$-generic 
and $q$ and $f$ are $P_\alpha$-names in $N$ 
and\/ {\rm ${\bf 1}\forces_{P_\alpha}``q\in P_{\alpha,\kappa}$ and 
$f$ is a $P_{\alpha,\kappa}$-name
and $q\forces_{P_{\alpha,\kappa}}`f\in{}^\omega\omega\cap V[G_{P_\alpha}]$,'\thinspace''}  
then there are $q'$ and $H$ such that $q'$ is a $P_\alpha$-name and
$H$ is a $P_\kappa$-name and\/ {\rm
$p\forces``q'\leq q$ and 
$ q'
\forces`H\in{}^\omega(\omega^{<\omega})\cap V[G_{P_{\alpha}}]$ and 
 $(\forall n)\allowbreak\allowbreak(f(n)\in H(n))$.'\thinspace''}

Proof: The proof proceeds by induction on $\kappa$.  We assume that 
  $\lambda$, 
$N$, $\alpha$, $p$, $q$, and $f$ are as in the hypothesis of the Theorem.

Case 1. $\kappa=\beta+1$.

Using the fact that ${\bf 1}\forces_{P_\beta}``Q_\beta$ has the Sacks
property,'' take $\tilde q$ and $H'$ to be $P_\beta$-names such that 
${\bf 1}\forces_{P_\beta}``\tilde q\leq q(\beta)$ and $H'$ is a $Q_\beta$-name and
$\tilde q\forces``H'\in{}^\omega(\omega^{<\omega})\cap V[G_{P_\beta}]$ and
$(\forall n\in\omega)\allowbreak(f(n)\in H'(n))$.'\thinspace''
We may assume that the names $\tilde q$ and $H'$ are elements of $N$.
Using Lemma 3.2 we may take $H^*$ and $q^*$ to be $P_\alpha$-names such that
${\bf 1}\forces_{P_\alpha}``H^*\in{}^{\omega}(\omega^{<\omega})$ and $q^*\leq q\restr\beta$ and
$q^*\forces_{P_{\alpha,\beta}}`H'=H^*$.'\thinspace'' We may assume that the names
$H^*$ and $q^*$ are in $N$.
Use the induction hypothesis to take a $P_\alpha$-name $q^\dag$ and
a $P_\beta$-name $G$  such that
$p\forces``q^\dag\leq q^*$ and 
$q^\dag\forces`G\in{}^\omega((\omega^{<\omega})^{<\omega})\cap V[G_{P_\alpha}]$ and
$(\forall n\in\omega)\allowbreak(H^*(n)\in G(n))$.'\thinspace''
 Let $H$ be 
a $P_\beta$-name such that
${\bf 1}\forces_{P_\beta}``H\in{}^\omega(\omega^{<\omega})$ and
$(\forall n\in\omega)\allowbreak(H(n)=\bigcup G(n))$.'' 
Let $q'$ be a $P_\alpha$-name such that ${\bf 1}\forces_{P_\alpha}``q'\in P_{\alpha,\kappa}$ and $q'\restr\beta=q^\dag$ and
$q^\dag\forces`q'(\beta)=\tilde q$.'\thinspace''
We have that $q'$ and $H$ are as required.

Case 2.  ${\rm cf}(\kappa)>\omega$.

Because no $\omega$-sequences of ordinals can be added at limit stages of uncountable cofinality, we may take $\beta$ and $f'$ and $q'$ to be $P_\alpha$-names in $N$ such that
${\bf 1}\forces``\alpha\leq\beta<\kappa$ and $f'$ is a
$P_{\alpha,\beta}$-name and $q'\leq q$ and
$q'\restr\beta\forces_{P_{\alpha,\beta}}`f'\in {}^\omega\omega$ and
$q'\restr[\beta,\kappa)\forces_{P_{\beta,\kappa}}``f'=f$.''\thinspace'\thinspace''

For every $\beta_0\in\kappa\cap N$ such that $\alpha\leq\beta_0$ 
let $\tilde q(\beta_0)$ and $H(\beta_0)$  be  $P_\alpha$-names in $N$ such that
${\bf 1}\forces``$if $\beta=\beta_0$ and there is some $\tilde q\leq q'\restr\beta$ and some 
$H^*$ such that $H^*$ is a $P_{\alpha,\beta}$-name and
$\tilde q\forces`H^*\in{}^\omega(\omega^{<\omega})\cap V[G_{P_\alpha}]$ and 
$(\forall n\in\omega)\allowbreak(f'(n)\in H^*(n))$,'
then $q^*(\beta_0)$ and $H(\beta_0)$  are witnesses thereto.''
Let $ q^*$ and $H$  and $s$ be $P_\alpha$-names such that for every $\beta_0
\in \kappa\cap N$, if $\alpha\leq\beta_0$, then
${\bf 1}\forces``$if $\beta=\beta_0$ then $ q^*=q^*(\beta_0)$ and $H=H(\beta_0)$
 and 
$s\in P_{\alpha,\kappa}$ and $s\restr\beta=q^*$ and
$s\restr[\beta,\kappa)=q'\restr[\beta,\kappa)$.''

Claim 1: $p\forces``s\leq q$ and $s\in N[G_{P_\alpha}]$ and $s
\forces`H\in{}^\omega(\omega^{<\omega})\cap V[G_{P_\alpha}]$ and
$(\forall n\in\omega)\allowbreak(f(n)\in H(n))$.'\thinspace''

Proof:  Suppose $p'\leq p$. Take $p^*\leq p'$ and $\beta_0<\kappa$ such that
$p^*\forces``\beta_0=\beta$.''  Because the name $\beta$ is in $N$ and $p^*$ is $N$-generic,
we have that $\beta_0\in N$.  Notice by the induction hypothesis
we have $p\forces``$there is some $q^\#\leq q'\restr\beta_0$ and 
some $P_{\alpha,\eta_0}$-names
$H^\#$ and $x^\#$
such that $q^\#\forces`H^\#\in{}^\omega(\omega^{<\omega})\cap V[G_{P_\alpha}]$ and
$(\forall n\in\omega)\allowbreak(f'(n)\in H^\#(n))$.'\thinspace''
Hence $p^*\forces``q^*=q^*(\beta_0)\leq q'\restr\beta$ and $H=H(\beta_0)$ and
$q^*\forces`H\in{}^\omega(\omega^{<\omega})\cap V[G_{P_\alpha}]$ and $(\forall
n\in\omega)\allowbreak(f'(n)\in H(n))$ and
$q'\restr[\beta,\kappa)\forces``f'=f$.''\thinspace'\thinspace''  Therefore
$p^*\forces``s\forces`(\forall n\in\omega)\allowbreak(f'(n)\in H(n))$.'\thinspace''

Claim 1 is established.
This completes Case 2.

Case 3. ${\rm cf}(\kappa)=\omega$.

Let $\langle\alpha_n\,\colon\allowbreak n\in\omega\rangle$ be an increasing 
sequence from $\kappa\cap N$ cofinal in $\kappa$ such that 
$\alpha_0=\alpha$.

Fix $\langle (p_n,f_n)\,\colon\allowbreak n\in\omega\rangle\in N$ (that is, the sequence of names is an element of $N$ but not necessarily their values) as in
Lemma 3.1 (applied in $V[G_{P_\alpha}])$.  That is, ${\bf 1}\forces``p_0\leq q$'' and 
 for every $n\in\omega$ we have that each of the following holds:

(0) $p_n$ is a $P_\alpha$-name for an element of $P_{\alpha,\kappa}$.

(1) For every $k\leq n$ we have ${\bf 1}\forces_{P_{\alpha_n}}``p_0\restr[\alpha_n,\kappa)
\forces`f(k)=f_n(k)$.'\thinspace''

(2) $f_n$ is a $P_{\alpha_n}$-name for an element of ${}^\omega\omega$.

(3) ${\bf 1}\forces_{P_{\alpha_n}}``p_0\restr[\alpha_n,\alpha_{n+1})\forces`f_n(k)=f_{n+1}(k)$ for every $k\leq n+1$.'\thinspace''

(4) ${\bf 1}\forces_{P_\alpha}``p_{n+1}\leq p_n$.''

(5) Whenever $k\leq m<\omega$ we have ${\bf 1}\forces_{P_{\alpha_n}}``p_m\restr[
\alpha_n,\alpha_{n+1})\forces`f_n(k)=f_{n+1}(k)$.'\thinspace''

Claim 2.  Suppose $\alpha\leq\beta\leq\gamma<\kappa$ and $\{\beta,\gamma\}\in N$, and
suppose $f'\in N$ is a $P_\gamma$-name for an element of ${}^\omega\omega$. Suppose
$r\in P_\beta$ is $N$-generic and $r\forces``q'\in P_{\beta,\gamma}\cap N[G_{P_\alpha}]$.''
Then there are $q^*$ and $H$ such that $r\forces``q^*\leq q'$ and
$q^*\forces`H\in {}^\omega(\omega^{<\omega})\cap V[G_{P_\alpha}]$ and $
(\forall i\in\omega)\allowbreak(f'(i)\in H(i))$.'\thinspace''

Proof: In $V[G_{P_\alpha}]$ let $D=\{s\in P_{\alpha,\beta}\,\colon\allowbreak
s\forces``(\exists q^*\leq q')\allowbreak(\exists H\in{}^\omega(\omega^{<\omega})\cap 
V[G_{P_\alpha}])
\allowbreak(q^*\forces_{P_{\beta,\gamma}}`(\forall 
i\in\omega)\allowbreak(f'(i)\in H(i))$')''$\}$.

Subclaim 1.  In $V[G_{P_\alpha}]$ we have that $D$ is a dense subset of
$P_{\alpha,\beta}$.

Proof: Given $\tilde p\in P_\alpha$ and $s$ a $P_\alpha$-name such that $\tilde
p\forces``s\in P_{\alpha,\beta}$,''
take $\lambda^*$ a sufficiently large regular cardinal and $M$ a countable elementary substructure of $\lambda^*$ containing $\tilde p$, $s$, and $N$.  Choose $p'\leq 
\tilde p$ such that $p'$ is $M$-generic.  By the overall induction hypothesis 
(i.e., because $\gamma<\kappa$), we have
$p'\forces``(\exists s'\in P_{\alpha,\gamma})\allowbreak(\exists H)\allowbreak
(s'\leq (s,q')$ and $H$ is a $P_{\alpha,\gamma}$-name and $s'\forces`H\in{}^\omega
(\omega^{<\omega})
\cap V[G_{P_\alpha}]$ and $(\forall i\in\omega)\allowbreak(f'(i)\in H(i))$').''  
Fix such $s'$ and $H$.
By Lemma 3.2 we may take $H'$  and $\tilde s$ to be $P_\alpha$-names such that
$p'\forces``\tilde s\leq s'$ and 
$h'\in{}^\omega\omega$ and $\tilde s\forces`H'=H$.'\thinspace''
We have $p'\forces``\tilde s\restr[\alpha,\beta)\in D$.''
The Subclaim is established.

Working in $V[G_{P_\alpha}]$, let ${\cal J}\subseteq D$ be a maximal antichain 
of $P_{\alpha,\beta}$.
For each $s\in{\cal J}$ take $q^*(s)$ and $H(s)$ witnessing that $s\in D$.
Construct $q^\#$ and $H^\#$ such that $(\forall s\in{\cal J})\allowbreak
(s\forces_{P_{\alpha,\beta}}``q^\#=q^*(s)$ and $H^\#=H(s)$''). 

Subclaim 2.  $r\forces``q^\#\leq q'$ and $q^\#\forces_{P_{\beta,\gamma}}
`(\forall i\in\omega)\allowbreak(f'(i)\in H^\#(i))$.'\thinspace''

Proof: Given $r_1\leq r$, take $r_2\leq r_1$ and a $P_\alpha$-name $s$ such that
$r_2\restr\alpha\forces``r_2\restr[\alpha,\beta)\leq s$ and $s\in{\cal J}$.''
We have $r_2\forces``q^\#=q^*(s)\leq q'$ and
$q^\#\forces`(\forall i\in\omega)\allowbreak
(f'(i)\in H(s)(i)=H^\#(i))$.'\thinspace''

The Subclaim and the Claim are established.

In $V[G_{P_\alpha}]$, fix $y$ and $z$ elements of
${}^\omega(\omega-\{0\})$ both eventually dominating every element of
${}^\omega(\omega-\{0\})\cap N[G_{P_\alpha}]$ such that $y\ll z$.
Using Lemma 4.8, build 
$\langle x^*_n\,\colon\allowbreak n\in\omega\rangle$ and
$\langle x_n\,\colon\allowbreak n\in\omega\rangle$
 sequences of
elements of ${}^\omega(\omega-\{0\})\cap N[G_{P_\alpha}]$ 
such that for every $n\in\omega$ we have
$x^*_n\ll x_n\ll x^*_{n+1}\ll y$ and $x^*_n$ eventually dominates every
element of
${}^\omega(\omega-\{0\})\cap N[G_{P_\alpha}]$.

For each $n\in\omega$ let $\langle T_{n,j}\,\colon\allowbreak j\in\omega\rangle$ list all
$x^*_{n+1}$-sized trees in $N[G_{P_\alpha}]$, and build $\langle x_{n,j}\,\colon\allowbreak j\in\omega\rangle$ a sequence of elements of
${}^\omega(\omega-\{0\})$ such that for every $j\in\omega$ we have that
$x_n\ll x_{n,j}\ll x_{n,j+1}\ll x^*_{n+1}$ and $T_{n,j}$ is an $x_{n,j}$-sized tree.

Using Lemma 4.7, take $T^*$ a $z$-sized tree and
 $\langle T^n\,\colon\allowbreak n\in\omega\rangle$ a sequence of trees such that
$T^0\subseteq T^*$ and $f_0\in[T^0]$ and for every $n\in\omega$ we have that
$T^n$ is an $x_n$-sized tree and $T^n\subseteq T^{n+1}$ and for every $j\in\omega$ and every 
$g\in[T_{n,j}]$  there is $k\in\omega$ such that for
every $\eta\in{}^{<\omega}\omega$ extending $g\restr k$, if
$\eta\in T_{n,j}$ and $\eta\restr k\in T^n\cap T^*$ then
$\eta\in T^{n+1}\cap T^*$.

Claim 3.  We may be build $\langle r_n\,\colon\allowbreak n\in\omega\rangle$ such that
$r_0=p$  and for every $n\in\omega$ we have that the following hold:

(1) $r_n\in P_{\alpha_n}$ is $N$-generic.

(2) $r_{n+1}\restr\alpha_n=r_n$.

(3) $r_n\forces``f_n\in[T^n]\cap[T^*]$.''

(4) $p\forces``r_n\restr[\alpha,\alpha_n)\leq p_0\restr\alpha_n$.''

Proof: Suppose we have $r_n$.

Let $F_0$ and $F_2$ be $P_{\alpha_n}$-names  such that
 ${\bf 1}\forces``$if there are functions $F'_0$ and
$F'_2$ such that $(\forall q'\in P_{\alpha_n,\alpha_{n+1}}
 )\allowbreak (F'_0(q')\in{}^\omega(\omega^{<\omega})\cap V[G_{P_\alpha}]$
 and
$F'_2(q')\leq q'$ and 
$F'_2(q')\forces`(\forall i\in\omega)\allowbreak(f_{n+1}(i)\in F'_0(q')(i))$'),
then $F_0$ and $F_2$ are witnesses.'' We may assume that the names $F_0$ and
$F_2$ are in $N$.

By Claim 2, we have that
$r_n\forces``(\forall q'\in P_{\alpha_n,\alpha_{n+1}}
 )\allowbreak (F_0(q')\in{}^\omega(\omega^{<\omega})\cap V[G_{P_\alpha}]$
 and
$F_2(q')\leq q'$ and 
$F_2(q')\forces`(\forall i\in\omega)\allowbreak(f_{n+1}(i)\in F_0(q')(i))$').''

For each $i\in\omega$ define $\tilde T(i)$ and $x_i$ and $T(i)$ 
as follows. In $V[G_{P_{\alpha_n}}]$, 
define $\tilde T(i)\subseteq {}^{<\omega}\omega$ by
$\eta\in\tilde T(i)$ iff $(\forall t<{\rm length}(\eta))\allowbreak
(\eta(t)\in F_0(p_i\restr[\alpha_n,\alpha_{n+1})))$.

Subclaim 1. For every $i\in\omega$ we have (in $V[G_{P_{\alpha_n}}]$)
that $f_n\restr i\in \tilde T(i)$.

Proof: Work in $V[G_{P_{\alpha_n}}]$. Because $F_2(p_{i}\restr[\alpha_n,\alpha_{n+1}))\leq
p_{i}\restr[\alpha_n,\alpha_{n+1})$ and $p_{i}
\restr[\alpha_n,\alpha_{n+1})\forces``f_n\restr i
=f_{n+1}\restr i$,'' we have that $F_2(p_{i}\restr[\alpha_n,\alpha_{n+1}))\forces``
f_n\restr i\in\tilde T(i)$.'' Therefore outright in $V[G_{P_{\alpha_n}}]$ we have
that $f_n\restr i\in\tilde T(i)$.  The Subclaim is established.

Working in $N[G_{P_{\alpha_n}}]$,
for each $i\in\omega$, take $x_i'\in{}^\omega(\omega-\{0\})$ such that
$\tilde T(i)$ is an $x'_i$-sized tree.  Define $y'\in{}^\omega(\omega-\{0\})$ by
$y(i)={\rm max}\{x'_j(i)\,\colon j\leq i\}$.
Let $T(i)=\{\eta\in\tilde T(i)\,\colon\allowbreak\eta\restr i= f_n\restr i\}$.
  Take
$z'\in{}^\omega(\omega-\{0\})$ such that $y'\ll z'$.

Working in $N[G_{P_{\alpha_n}}]$, we know that
 for every $i\in\omega$ we have that
$T(i)$ is an $y'$-sized tree and $h_i\in[T(i)]$ and $f_n\restr i=h_i\restr i$.
Therefore we may use Lemma 4.6 to take $\tilde T$ a $z'$-sized tree and
$\langle k_i\,\colon\allowbreak i\in\omega\rangle$ an increasing sequence of
integers such that for every $\eta\in{}^{<\omega}\omega$ and every $i\in\omega$ and every
$\nu\in{}^{<\omega}\omega$, if $\eta$ is an
initial segment of $f_n$ and ${\rm length}(\eta)\geq k_i$ and $\nu\in T(k_i)$ and
$\nu$ extends $\eta$, then $\nu\in\tilde T$.

Subclaim 2.  For all $i\in\omega$ we have that

\medskip

{\centerline{$r_n\forces``F_2(p_{k_i}\restr[\alpha_n,\alpha_{n+1}))
\forces`f_{n+1}\in[\tilde T]$.'\thinspace''}}

\medskip

Proof. Working in $V[G]$ where
$G$ is $P_{\alpha_{n+1}}$-generic and
$(r_n,F_2(p_{k_i}\restr[\alpha_n,\alpha_{n+1})))\in G$, we may
 let $\eta=f_n\restr k_i$ and
let $\nu$ be any initial segment of $f_{n+1}$ extending $\eta$.
By the definition of $F_2$ we have $f_{n+1}\in[\tilde T(k_i)]$, and therefore
$f_{n+1}\in[T(k_i)]$, and therefor $\nu\in T(k_i)$.
Looking at the definition of $\tilde T$ we see that $\nu\in\tilde T$. The
Subclaim is established.

Subclaim 3.  ${\bf 1}\forces_{P_{\alpha_n}}``f_n\in[\tilde T]$.''

Proof.  Working in $V[G_{P_{\alpha_n}}]$, for each 
$i\in\omega$, let $\eta_i=\nu_i=f_n\restr k_i$.  By the definition of
$\tilde T$ (with $\eta=\eta_i$ and $\nu=\nu_i$)
we see that $\eta_i\in \tilde T$.  The Subclaim is established.

Notice that $N[G_{P_{\alpha_n}}]\models``\tilde T\in V[G_{P_\alpha}]$,'' hence using 
Lemma 3.3 we have
$r_n\forces``\tilde T\in N[G_{P_\alpha}]$.''

Let $k^*$ be a $P_{\alpha_n}$-name for an integer such that
in $V[G_{P_{\alpha_n}}]$ we have $$(\forall t\geq k^*)(\vert \tilde T\cap{}^t\omega\vert
\leq x^*_{n+1}(t))$$
In $N[G_{P_{\alpha_n}}]$ define $T^\dag\subseteq {}^{<\omega}\omega$ by
$T^\dag=\{\eta\in\tilde T\,\colon\allowbreak\eta\restr k^*=f_n\restr k^*\}$.

We have $r_n\forces``T^\dag\in N[G_{P_\alpha}]$'' and $T^\dag$ is an $x^*_{n+1}$-sized tree,
so
 we may choose a $P_{\alpha_n}$-name $k$ such that $r_n\forces`` T^\dag
=T_{n,k}$.''

Subclaim 4.  ${\bf 1}\forces_{P_{\alpha_n}}``f_n\in [T^\dag]$.''

Proof:  Immediate from Subclaim 3 and the definition of $T^\dag$.

Using Subclaim 4, the fact that 
Claim 3 holds for the integer $n$, the fact that $r_n\forces``T^\dag=T_{n,k}$,'' 
and the fact that the trees $\langle T^n\,\colon\allowbreak n\in\omega\rangle$ and
$T^*$ were chosen as in the hypothesis of Lemma 4.7,
 we may, by Lemma 4.7,  choose
$k'$ a $P_{\alpha_n}$-name for an integer such that 

\medskip

(*) $r_n\forces``(\forall\eta\in T^\dag)\allowbreak($if $\eta$ extends $f_n\restr k'$ and
$\eta\restr k'\in T^n\cap T^*$ then $\eta\in T^{n+1}\cap T^*)$.

\medskip

Choose $j\in\omega$ such that $k_j\geq{\rm max}(k', k^*)$. Let $K=k_j$.

Subclaim 5. 
$r_n\forces``F_2(p_K\restr[\alpha_n,\alpha_{n+1}))\forces
`f_{n+1}\in [T^{n+1}]\cap[T^*]$.'\thinspace''

Proof.  By Subclaim 2 we have 
$r_n\forces``F_2(p_K\restr[\alpha_n,\alpha_{n+1}))\forces
`f_{n+1}\in [\tilde T]$.'\thinspace''
Because $r_n\forces``F_2(p_K\restr[\alpha_n,\alpha_{n+1}))\leq 
p_{k^*}\restr[\alpha_n,\alpha_{n+1})$,'' we have

\medskip

(**) $r_2\forces``F_2(p_K\restr[\alpha_n,\alpha_{n+1}))\forces`f_{n+1}\in[T^\dag]$.'\thinspace''

\medskip

By (*) and (**) we have 

\medskip

$r_n\forces``F_2(p_K\restr[\alpha_n,\alpha_{n+1}))\forces`(\forall i\in\omega)\allowbreak($if
$f_{n+1}\restr i$ extends $f_n\restr k'$ and $f_{n+1}\restr k'\in T^n\cap T^*$, then
$f_{n+1}\restr i\in T^{n+1}\cap T^*$.'\thinspace''

\medskip

Therefore, to establsih the Subclaim, it suffices to show

\medskip

{\centerline{$r_n\forces``F_2(p_K\restr[\alpha_n,\alpha_{n+1}))
\forces`f_{n+1}\restr k'=f_n\restr k'$' and
$f_n\restr k'\in T^n\cap T^*$.''}}

\medskip

The first conjunct follows from the fact that
 $r_2\forces``F_2(p_K\restr[\alpha_n,\alpha_{n+1}))\leq
p_{k'}\restr[\alpha_n,\alpha_{n+1})$,'' and the second conjunct follows from the
assumption that Claim 3 holds for all integers less than or equal to $n$.
The Subclaim is established.

Use the Proper Iteration Lemma to take $r_{n+1}\in P_{\alpha_{n+1}}$ such that
$r_{n+1}\restr\alpha_n=r_n$ and $r_{n+1}$ is $N$-generic and
$r_n\forces``r_{n+1}\restr[\alpha_n,\alpha_{n+1})\leq F_2(p_K\restr[\alpha_n,\alpha_{n+1}))$.''

This completes the induction proving Claim 3.

Let $q'$  and $H$ be $P_\alpha$-names such that
$p\forces``q'=\bigcup\{r_n\restr[\alpha,\alpha_n)\,\colon\allowbreak n\in\omega\}$, and
$H\in{}^\omega(\omega^{<\omega})$ and
for all $i\in\omega$ we have
$H(i)=\{\eta(i)\,\colon\allowbreak \eta\in T^*$ and $i\in {\rm domain}(\eta)\}$.''
By Claim 3, we have that $q'$ and $H$ satisfy the requirements of the Theorem.

The Theorem is established.

\proclaim Corollary 4.10.  Suppose $\langle P_\eta\,\colon\allowbreak\eta\leq\kappa\rangle$ is a countable support iteration based on $\langle Q_\eta\,\colon\allowbreak\eta<\kappa\rangle$ and suppose\/ {\rm $(\forall\eta<\kappa)\allowbreak({\bf 1}\forces_{P_\eta}``Q_\eta$ 
is proper and has the Sacks proeprty.'')}  Then $P_\kappa$ has the Sacks property.

Proof.  Take $\alpha = 0$ in Theorem 4.9.

\end{document}